\documentclass[11pt, a4paper]{amsart}

\usepackage[english]{babel}
\usepackage[utf8]{inputenc} 

 \usepackage[T1]{fontenc}
 \usepackage{lmodern}

\usepackage{amsmath}
\usepackage{amsfonts}
\usepackage{amssymb}

\usepackage{amsthm}
\theoremstyle{plain}
\newtheorem{thm}{Theorem}
\newtheorem{prop}[thm]{Proposition}
\newtheorem{lemme}[thm]{Lemma}

\theoremstyle{definition}
\newtheorem{defn}[thm]{Definition}

\title{Corrugation Process and $\epsilon$-Isometric Maps}
\author{Mélanie \textsc{Theillière}}
\address{Institut Camille Jordan, Braconnier\\
Université Claude Bernard, Lyon 1\\
43 boulevard du 11 novembre 1918\\
F-69622 Villeurbanne Cedex}
\email{melanie.theilliere@univ-lyon1.fr}


\usepackage{xcolor}

\newcommand{\bg}{\raisebox{1.7pt}{$\gamma$\hspace{-4.5pt}\raisebox{-6.5pt}{\textcolor{white}{$\blacksquare$}}} \hspace{-18.25pt} $\gamma$}

\def\bgamma{\textstyle{\textup{\bg}}}
\def\hgamma{\textstyle{\textup{$\gamma$}}}

\usepackage[cal=boondox]{mathalfa}
\usepackage{graphicx}

\def\R{\mathbb{R}}
\def\Z{\mathbb{Z}}
\def\C{\mathbb{C}}

\def\E{\mathbb{E}}

\def\proj{\mbox{proj}}

\def\rel{\mathcal{R}}
\def\is{\mathcal{I \!\! s}}
\def\dj{\partial_j}

\def\noi{\noindent}

\newcommand{\CP}{CP_{\hgamma}}
\newcommand{\s}{\mathfrak{S}}
\newcommand{\iconv}{\mathrm{IntConv}}

\begin{document}

\begin{abstract}
Convex Integration is a theory developed in the '70s by M. Gromov. This theory allows to solve families of differential problems satisfying some convex assumptions. From a subsolution, the theory iteratively builds a solution by applying a series of \textit{convex integrations}. In a previous paper \cite{the01}, we proposed to replace the usual convex integration formula by a new one called \textit{Corrugation Process}. This new formula is of particular interest when the differential problem under consideration has the property of being \textit{of Kuiper}. In this paper, we consider the differential problem of $\epsilon$-isometric maps and we prove that it is Kuiper in codimension 1. As an application, we construct $\epsilon$-isometric maps from a short map having a conical singularity.
\end{abstract}

\maketitle

\section{General introduction}

\subsection{The Nash-Kuiper Theorem}

A map $f:(M,g)\rightarrow \E^n$ between a Riemannian manifolds $(M,g)$ and the Euclidean space $\E^n=(\R^n,\langle\cdot,\cdot\rangle)$ is said to be isometric if $g=f^*\langle\cdot,\cdot\rangle.$ It is said to be \textit{strictly short} if $g -f^*\langle\cdot,\cdot\rangle$ is positive definite (as usual $f^*h$ denotes the pullback of the metric $h$ by $f$). In other words, the length of the image of any curve in $M$ by a strictly short map is shorter than the length of the curve in $M$. The $C^1$ embedding theorem of Nash and Kuiper states that close to every strictly short map lies a $C^1$-isometric map:

\begin{thm}[\cite{Nash54,Kuiper55}]\label{thNK}
Let $(M^m,g)$ be a compact Riemannian manifold and let $f_0:(M^m,g)\rightarrow (\R^{n},\langle\cdot,\cdot\rangle)$, with $n>m$, be a strictly short embedding. Then for any $\epsilon >0$ there exists a $C^1$-isometric embedding $f:(M^m,g)\rightarrow (\R^n, \langle\cdot,\cdot\rangle)$ such that, $\sup_{x}\|f(x)-f_0(x)\| < \epsilon.$
\end{thm}

The proof considers an increasing sequence of metrics $g_k$ converging toward $g$ and a decreasing sequence $\epsilon_k$ converging toward $0$. A sequence of maps $f_1, \ldots$, $f_k$, $\ldots$ is then iteratively built such that, for each $k$, $f_k$ is an $\epsilon_k$-isometric map from $(M,g_k)$ to $\E^n$ i. e.
$$\|g_k- f_k^*h\| < \epsilon_k.$$
Parameters of the construction are chosen to insure the $C^1$ convergence of the sequence $(f_k)_k$ so that the limit map $f_{\infty}$ is isometric.

\subsection{Differential relations}\label{subsect_diff_rel}

We now introduce the formalism of Gromov's Convex Integration Theory~\cite{Gromov-book}. This theory can be seen as a wide generalization of Nash's approach. It provides a powerful tool to solve a large family of differential constraints. We denote by
\begin{eqnarray*}
	J^1(M,W) := \{(x,y,L) \,\, | \,\, x\in M,\,\, y\in W,\,\, L:T_xM\rightarrow T_yW \mbox{ a linear map} \}.
\end{eqnarray*}
the $1$-jet space of $C^1$ maps between $M$ and $W$.  Every $C^1$-map $f$ gives rise to a section $j^1f:x\mapsto (x,f(x), df_x)$ of $J^1(M,W)$ called the \textit{1-jet} of $f$. For any section $\s_0 :x\mapsto (x, f_0(x),L_x)$ we denote by $f_0= \mathrm{bs}\;\; \s_0$ its \textit{base map}.\\

A \textit{differential relation} $\rel$ is any subset of the $1$-jet space $J^1(M,W)$. For instance, the $\epsilon$-isometric condition defines the differential relation $\is(\epsilon)$ of $\epsilon$-isometric maps:
\begin{eqnarray*}
\is(\epsilon) := \{ (x,y,L) \,\, | \,\,  \| g_x - L^*h_y \| < \epsilon \}
\end{eqnarray*}
where $L^*h$ denotes the pullback by $L$ of the metric $h$. Observe that, as a topological subspace, $\is(\epsilon)$ is open. 

\begin{defn}
Let $\s:M \rightarrow J^1(M,W)$ be a section. We say that $\s$ is a \textit{formal solution of $\rel$} if the image of $\s$ lies in~$\rel$. Moreover, if there exists a $C^1$-map $f:M \rightarrow W$ such that $j^1 f = \s$, we say that $\s$ is a \textit{holonomic solution of $\rel$}.
\end{defn}

For instance, building an $\epsilon$-isometric map $f$ is equivalent to finding a holonomic solution $j^1f$ of $\is(\epsilon)$.\\

Under some topological and convex assumptions on $\rel$, the Convex Integration Theory allows to deform a formal solution $\s_0$ to a holonomic one. Each convex integration modifies the 1-jet of a formal solution $(x_0,f_0,L_0)$ in a given direction $u$ to obtain a new formal solution $\s=(x,f,L)$ such that $L(u)=df(u)$. Loosely speaking $\s$ is "partially holonome in the direction $u$". Here is how it works for $M=[0,1]^m$ and $W=\E^n$. In this case a formal solution writes
$$\s_0 : x \mapsto (x, f_0(x), v_1(x),\ldots, v_m (x))\in \rel$$
where we have identified the 1-jet space with the product
$$J^1([0,1]^m,\R^n)=[0,1]^m\times\R^n\times (\R^n)^m$$
and $\s_0$ is holonomic if there exists a map $f$ such that
$$\s_0 = j^1 f :x \mapsto (x, f(x), \partial_1 f(x),\ldots, \partial_m f(x)).$$
From a formal solution $\s_0$, the Convex Integration Theory builds a finite sequence of formal solutions $\s_k$ such that for every $k\in\{1,\ldots,m\}$ we have
\begin{eqnarray*}
\s_k : x \mapsto (x, f_k(x), \partial_1 f_k(x),\ldots, \partial_k f_k(x), v_{k+1}(x),\ldots, v_m(x)) \in \rel.
\end{eqnarray*}
In particular
$$\s_m=j^1 f_m : x \mapsto (x, f_m(x), \partial_1 f_m(x),\ldots,\partial_m f_m(x))$$
is a holonomic solution of $\rel$.

\subsection{Corrugation Process}

To build the sequence $(\s_k)_k$, we propose in~\cite{the01} to replace the usual formula of the Convex Integration Theory by another one, called Corrugation Process:

\begin{defn}\label{def_formula_CP}
Let $f_0:[0,1]^m\rightarrow \mathbb{R}^n$ be a map, $\dj$ be a direction, $\gamma :[0,1]^m\times\R/\Z \rightarrow \mathbb{R}^n$ be a loop family and $N\in\; ]0,+\infty[$.
We define the map $f_1:[0,1]^m\rightarrow \mathbb{R}^n$ by 
\begin{eqnarray}\label{formula_CP}
f_1(x):=f_0(x) + \frac{1}{N} \int_{s=0}^{Nx_j} \gamma(x,s)-\overline{\gamma}(x) ds
\end{eqnarray}
where $\overline{\gamma}(x)=\int_0^1 \gamma(x,s) ds$ denotes the average of the loop $t\mapsto\gamma(x,t)$. We say that $f_1$ is obtained from $f_0$ by a \textit{Corrugation Process} in the direction $\partial_j$ and we denote $f_1=CP_{\hgamma}(f_0,\dj,N)$.
\end{defn}

The Corrugation Process satisfies the following three properties which are at the basis of the Convex Integration Theory~\cite{Spring-book}.

\begin{prop}[\cite{the01}]\label{propCP}
The map $f_1=CP_{\hgamma}(f_0,\dj,N)$ satisfies
\begin{itemize}
\item[$(P_1)$] \ \ $\|f_0-f_1\|_{C^0} = O(1/N)$,
\item[$(P_2)$] \ \ $\| \partial_{i} f_0 -\partial_{i} f_1 \|_{C^0} = O(1/N)$ for every $i\neq j$.
\end{itemize}
Moreover if $\forall x\in [0,1]^m$ we have $\partial_{j} f_0(x) = \overline{\gamma}(x)$ then
\begin{itemize}
\item[$(P_3)$]\ \ $\partial_{j} f_1(x) = \gamma(x,Nx_j) + O(1/N)$ for all $x\in [0,1]^m$.
\end{itemize}
\end{prop}

Provided that $N$ is large enough, this proposition shows that the Corrugation Process allows to modify the $j$th partial derivative while keeping the other derivatives under control. Consequently, this formula performs the deformations required to build the sequence $(\s_k)_k$ provided $\gamma$ has values in a well chosen region. For more details and for a proof of this proposition see~\cite{the01}. We give below a coordinate free expression of the Corrugation Process:

\begin{defn}
Let $f_0:U\rightarrow (W,h)$ be a map from an open set $U\subset M$, $\pi:U\rightarrow \R$ be a submersion and $\gamma:U\times \R/\Z\rightarrow f_0^*TW$ be a loop family such that $\gamma(x,.):\R/\Z\rightarrow f_0^*TW_x$ for every $x\in U$. The map defined by Corrugation Process is defined by
\begin{eqnarray*}
f_1 = \CP(f_0,\pi,N):x\mapsto \exp_{f_0(x)} \Big( \frac{1}{N}\int_{t=0}^{N\pi(x)} \gamma(x,t)-\overline{\gamma}(x) dt  \Big)
\end{eqnarray*}
where $\exp: TW \rightarrow W$ is the exponential map induced by the metric $h$.
\end{defn}

\subsection{Subsolutions}\label{subsect_subsol}
Subsolutions are a refinement of the notion of formal solution. This refinement is needed to ensure the existence of a loop family $\gamma$ whose its values is chosen in an appropriate region and whose its average is the partial derivative $\partial_jf_0$ under consideration (see property $(P_3)$ of Proposition~\ref{propCP}).\\

Let $\rel$ be a differential relation, $\sigma=(x,y,L)\in \rel$ and $(\lambda, u)\in T_x^*M\times T_xM$ such that $\lambda(u)=1$. We set
$$\rel(\sigma,\lambda,u) := Conn_{L(u)}\{ v\in T_yW \; |\; (x,y,L+(v-L(u))\otimes\lambda)\in \rel\}$$
where $Conn_aA$ denotes the path connected component of $A$ that contains~$a$. We say that $\rel(\sigma,\lambda,u)$ is the \textit{slice} of $\rel$ over $\s$ with respect to $(\lambda,u)$. Note that the linear map $L+(v-L(u))\otimes\lambda$ coincides with $L$ over $\ker \lambda$ and maps $u$ to $v$. We then denote by $\iconv\;\rel(\sigma,\lambda,u)$ the interior of the convex hull of $\rel(\sigma,\lambda,u).$

\begin{defn}\label{def_sous-sol_generale}
Let $U\subset M$, $\pi:U\rightarrow \R$ be a submersion and $u:U\rightarrow TM$ be a vector field such that $d\pi_x(u_x)=1$. Let $x\mapsto\s(x)=(x, f_0(x),L(x))$ be a formal solution of $\rel$ over $U$. If for every $x$ in $U$ the base map $f_0=\mathrm{bs}\;\s$ satisfies
$$df_0(u_x)\in \iconv\;\rel(\s(x),d\pi_x,u_x)$$
then the formal solution $\s$ is called a \textit{subsolution of} $\rel$ \textit{with respect to} $(d\pi,u)$.
\end{defn}

In the case where $U=[0,1]^m$, $W=\R^n$, $\pi(x)=x_j$ and $u=\partial_j$, the condition of the definition means that $\partial_j f_0(x)$ lies in the interior of the convex hull of 
$$\rel(\sigma,dx_j,\partial_j)=Conn_{v_j}\{ t\in\R^n\; |\; (x,f_0(x),v_1,\ldots,v_{j-1},t,v_{j+1},\ldots,v_m)\in\rel\}.$$

From a subsolution $\s$ of $\rel$ with respect to $(d\pi,u)$ the Convex Integration Theory builds a map $f_1$ whose derivative along $u_x$ lies in the slice $\rel(\s(x),d\pi_x,u_x)$:

\begin{lemme}\label{gamma_dans_R}
Let $\rel$ be an open differential relation and let $\s$ be a subsolution of $\rel$ with respect to $(d\pi,u)$ and with base map $f_0= \mathrm{bs}\;\; \s$. Then there exists a loop family $\gamma$ such that for every $x\in U$ we have $ \overline{\gamma}(x)=df_0 (u_x) $ and for every  $(x,t)\in U\times\R/\Z$ the image of $\gamma$ lies in $\rel(\s(x),d\pi_x,u_x)$. If we set $f_1:=\CP(f_0,\pi,N)$ for this loop family $\gamma$, we have
$$\forall x\in U,\;\;\;\;\;df_1(u_x)\in \rel(\s(x),d\pi_x,u_x)$$
for $N$ large enough.
\end{lemme}

\noi
\textbf{Proof.--} The existence of $\gamma$ follows the Integral Representation Lemma of the Convex Integration Theory of Gromov (\cite[p169]{Gromov-book} or \cite[p29]{Spring-book}). The property on $df_1 (u_x)$ is a direct consequence of point $(P_3)$ of Proposition~\ref{propCP}.
 \textcolor{white}{,}\hfill $\square$

\subsection{Kuiper relations}\label{subsection_Krel}

In the usual approach, the family of loops $\gamma$ is constructed a posteriori once the subsolution $\mathfrak{S}$ given. However the construction of a holonomic solution often requires to repeat the Corrugation Process in several directions $\partial_j$ and consequently needs to re-build at each step the loop family $\gamma$ on a different subsolution at each time. In \cite{the01}, we propose to simplify this approach by constructing a bigger loop family $\bgamma$ that could be used indifferently regardless of the subsolution. This simplification leads to introduce the notion of surrounding loop family and then the notion of Kuiper relation.\\

Basically, a surrounding family is a family of loops lying inside $\rel$ which is double indexed by its base point $\sigma$ and its average $w$ and where $(\sigma,w)$ are allowed to vary in the largest possible space, that is, inside 
$$\iconv(\mathcal{R},d\pi,u):=\{(\sigma,w)\in p_y^*TW\; |\; w\in \iconv\;\rel(\sigma,d\pi_x,u_x)\}.$$
In that definition, $p_y^*TW$ is the bundle over $\rel$ induced by the projection $p_y:\rel\rightarrow W$, $\sigma=(x,y,L)\mapsto y.$ 

\begin{defn}\label{def_surrounding}
Let $\rel$ be a differential relation of $J^1(U,W)$. We say that a loop family
\begin{eqnarray*}
\begin{array}{lclc}
\bgamma : & \iconv(\rel,d\pi,u) & \longrightarrow & C^{0}(\R/\Z,TW)\\
& (\sigma,w) & \longmapsto & \bgamma (\sigma,w)(\cdot)
\end{array}
\end{eqnarray*}
is {\it surrounding with respect to} $(d\pi,u)$ if for every $(\sigma,w)$ we have
\begin{itemize}
\item[$(1)$]\ \ $t\mapsto \bgamma (\sigma,w)(t)$ is a loop in $\mathcal{R}(\sigma,d\pi_x,u_x)$,
\item[$(2)$]\ \ the average of $t\mapsto \bgamma (\sigma,w)(t)$ is $w$,
\item[$(3)$]\ \ there exists a continuous homotopy $H:\iconv(\mathcal{R},d\pi,u)\times [0,1]\rightarrow TW$ such that
$H(\sigma,w,0)=\bgamma (\sigma,w)(0)$, $H(\sigma,w,1)=L(u_x)$ and $H(\sigma,w,t)\in \mathcal{R}(\sigma,d\pi_x,u_x)$ for all $t\in\, [0,1].$
\end{itemize} 
\end{defn}

Note that point $(3)$ is a homotopic property needed to state a potential $h$-principle for $\rel$.\\

Then for any subsolution $\s=(x,f_0,L)$ we choose the loop family
$$ \gamma(x,t):=\bgamma(\s(x), df_0(u_x))(t)\in \rel(\sigma,d\pi_x,u_x)$$
for every $(x,t)\in U\times\R/\Z$, and we write $CP_{\bgamma}(\s,\pi,N) := CP_{\hgamma}(f_0,\pi,N).$\\

We would like to ensure that all loops $\bgamma(\sigma,w)$ share the same pattern.

\begin{defn}\label{def_c}
Let $p,q>0$ be two natural numbers and $A\subset \R^q$ be a parameter space. A family of $1$-periodic curves $c:A\times \R /\Z\rightarrow\R^p$ is said to be a \textit{pattern}.
\end{defn}

We denote by $E\rightarrow W$ the fiber bundle over $W$ with fiber $\mathcal{L}(\R^p,T_yW)=(\R^p)^*\otimes T_yW$ and we consider its pull back by the projection $q:\iconv(\rel,d\pi,u)\rightarrow W$, $(\sigma,w)\mapsto y$. A section $\mathcal{e}$ of $q^*E$ defines a family of linear maps $\mathcal{e}(\sigma,w):\R^p\rightarrow T_yW$.

\begin{defn}\label{def_R_Kuiper}
Let $c$ be a loop pattern. If there exist a surrounding loop family $\bgamma:  \iconv(\rel,d\pi,u)  \rightarrow  C^{0}(\R/\Z,TW)$ with respect to $(d\pi,u)$, a section $\mathcal{e}$ of $q^*E\rightarrow \iconv(\mathcal{R},d\pi,u)$ and a map
${\bf a}:\iconv(\mathcal{R},d\pi,u)\rightarrow A$ such that, for all $((\sigma,w),t)\in \iconv(\mathcal{R},d\pi,u)\times \R/\Z$,
\begin{eqnarray*}
\bgamma(\sigma,w)(t) = \displaystyle \mathcal{e}(\sigma,w)\circ c({\bf a}(\sigma,w),t)
\end{eqnarray*}
we then say that $\rel$ is a {\it Kuiper relation with respect to $(c,d\pi,u)$}.
\end{defn}

If $(c_1,\ldots,c_p)$ denote the components of $c$ in the standard basis of $\R^p$ and if $\textbf{e}_1, \ldots, \textbf{e}_p$ denote the image of this basis by $\mathcal{e}$, the above definition writes
\begin{eqnarray*}
\bgamma(\sigma,w)(t) = \displaystyle \sum_{i=1}^p c_i({\bf a}(\sigma,w),t)\, {\bf e}_i(\sigma,w).
\end{eqnarray*}
We denote the periodic primitive of the $c_i$'s by
$$C_i(a,t)=\int_{s=0}^t c_i(a,s)- \overline{c_i}(a)ds.$$

\begin{prop}\label{prop_noInt}
Let $c$ be a loop pattern, $\rel$ be an open Kuiper relation with respect to $(c,d\pi,u)$, $\s=(x,f_0,L_0)$ be a subsolution and $\bgamma$ be a $c$-shaped surrounding loop family. Then $f_1 = CP_{\bgamma}(\s,\pi,N)$ has the following analytic expression
\begin{eqnarray}
f_1(x) &=&  \exp_{f_0(x)}\left(\frac{1}{N} \sum_{i=1}^p C_i( a(x),N\pi(x) ) e_i(x)\right)
\end{eqnarray}
where $a(x) := \textbf{a}(\s(x), df_0(u_x))$, $e(x) := \textbf{e}(\s(x), df_0(u_x))$ and $x\in U$. Moreover, if $N$ is large enough, the section
$$x\mapsto \s_1:=(x,f_1,L_1 = L_0+(df_1(u_x)-L_0(u_x))\otimes d\pi)$$
is a formal solution of $\rel$. 
\end{prop}

 In the case where $U=[0,1]^m$, $W=\R^n$, $\pi(x)=x_j$ and $u=\partial_j$ the map $f_1=CP_{\bgamma}(\s,\dj,N)$ is given by
\begin{eqnarray*}
f_1(x) = f_0(x) + \frac{1}{N} \Big( \sum_{i=1}^p C_i(a(x),Nx_j) e_i(x) \Big).
\end{eqnarray*}

In \cite{the01} the reader will find a proof of the proposition as well as examples of Kuiper relations. In the next section, we prove that the relation of $\epsilon$-isometric maps is Kuiper in codimension one.

\section{The relation of $\epsilon$-isometric maps}\label{section_is}

In this article, we prove the following theorem:

\begin{thm}\label{th_is_Kuiper}
Let $M$ and $W$ be orientable Riemannian manifolds such that $\dim W = \dim M +1$. For every $\epsilon>0$, the relation $\is(\epsilon)$ is a Kuiper relation.
\end{thm}

The key point of the proof of this theorem is to build a loop family $\bgamma$ $c$-shaped for all couples $(\sigma,w)$ such that $\sigma$ belongs to $\is(\epsilon)$ and $w$ belongs to the convex hull of the slice $\is(\epsilon)(\sigma,\lambda,u)$, for some $\lambda$, $u$. To understand the slice $\is(\epsilon)(\sigma,\lambda,u)$ and its convex hull, we first present its geometric description and a description of its subsolutions. We then give a proof of Theorem~\ref{th_is_Kuiper}.

\subsection{Geometric description of the relation of isometric maps}\label{subsect_geom_desc}

The relation of $\epsilon$-isometric maps is a thickening of the relation of isometric maps
\begin{eqnarray*}
\is := \{ (x,y,L) \,\, | \,\,   g = L^*h  \} \subset J^1(M,W)
\end{eqnarray*}
where $g$ is a metric of $M$ and $L^*h$ is the pullback by $L$ of the metric $h$ of $W$. So in this paragraph we give a geometric description of the relation of isometric maps. Such a description can be found in \cite[p202]{Gromov-book} or \cite[p194]{Spring-book}. For the sake of completeness we recall this description here in the coordinate-free case and we give some extra details needed for our construction of a surrounding loop family of the relation of $\epsilon$-isometric maps.\\

\noi 
Let $\sigma = (x,y,L)\in \is$. Let $\lambda \in T_x^* M$ and $u\in T_xM$ such that $\lambda(u)=1$. For every $v\in T_yW$, we set $L_v := L+(v-L(u))\otimes \lambda$. We have
\begin{eqnarray*}
\is(\sigma,\lambda, u) 
&:=&  Conn_{L(u)} \{ v \in T_yW \;|\; (x,y,L_v) \in \is \}\\
&=&  Conn_{L(u)} \{ v \in T_yW \;|\;g_x = L_v^*h_y \}.
\end{eqnarray*}

\begin{figure}[!ht]
\centering
\includegraphics[scale=0.52]{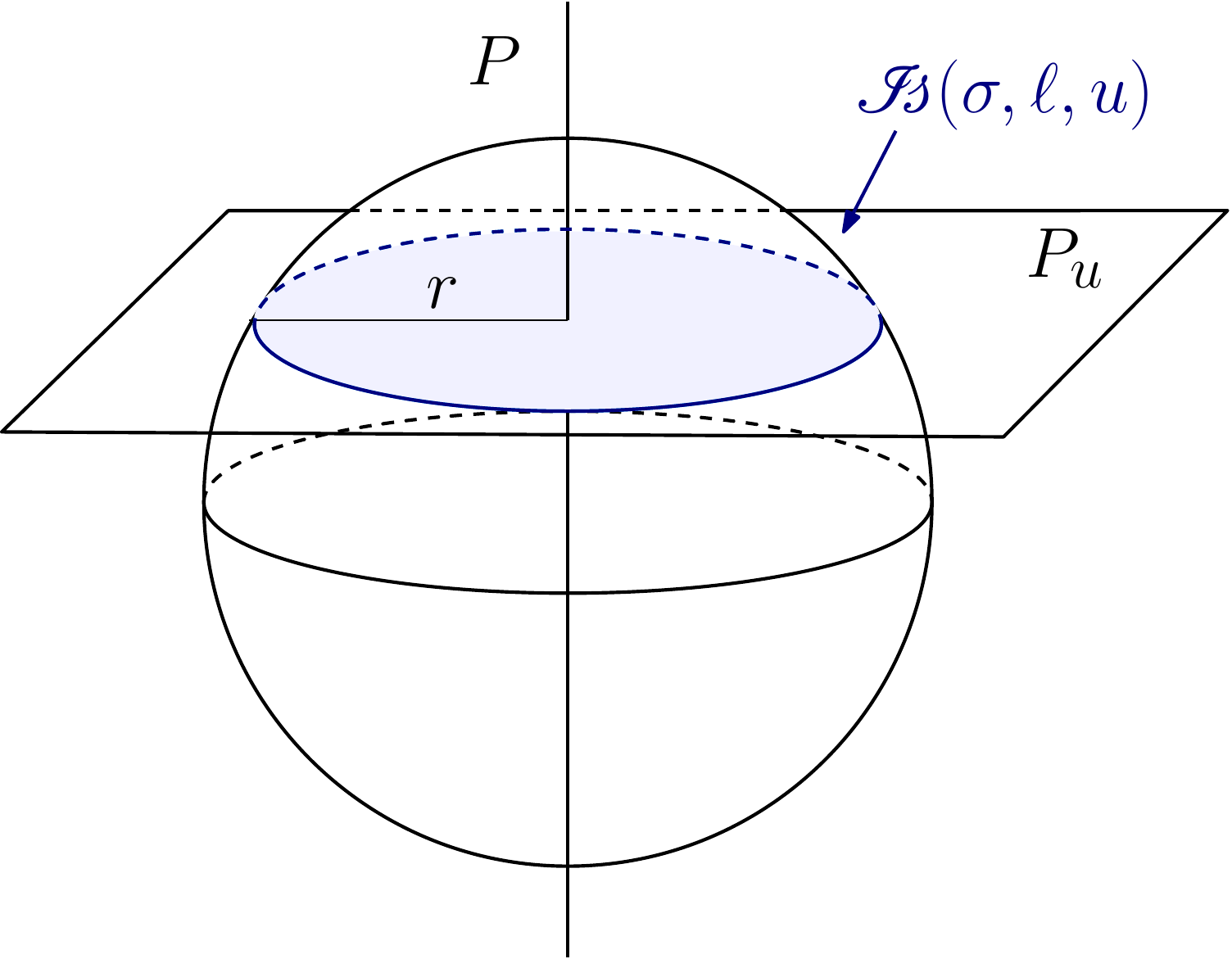}
\caption{\footnotesize {\textbf{The slice $\is(\sigma,\lambda,u)$ and its convex hull:} the $(n-m)$-dimensional sphere in dark blue is $\is(\sigma,\lambda,u)$ and the convex hull $\iconv \, \is(\sigma,\lambda,u)$ is the $(n-m+1)$-ball in light blue. $P$ denotes the $(m-1)$-plane $L(\ker \lambda)$}} \label{image_relation_isometrie}
\end{figure}

Note that, by the definition of $L_v$, we have $L_v(u) =v$ and for every $u_0\in \ker \lambda$ we have $L_v(u_0) = L(u_0)$, in particular $L_v(\ker \lambda) = L(\ker \lambda)$. Let $w_1= \alpha_1 u+ a_1$ and $w_2 = \alpha_2 u +a_2$ with $\alpha_1,\alpha_2\in \R$ and $a_1,a_2\in \ker \lambda$. As $g=L^*h$, we have
\begin{eqnarray*}
(g - L_v^*h)(w_1,w_2) &=& \alpha_1 \alpha_2 (g(u,u)-h(v,v))\\
&& + \alpha_1 h(L(u)-v,L(a_2)) + \alpha_2 h(L(u)-v, L(a_1)).
\end{eqnarray*}

\begin{figure}[!ht]
\centering
\includegraphics[scale=0.54]{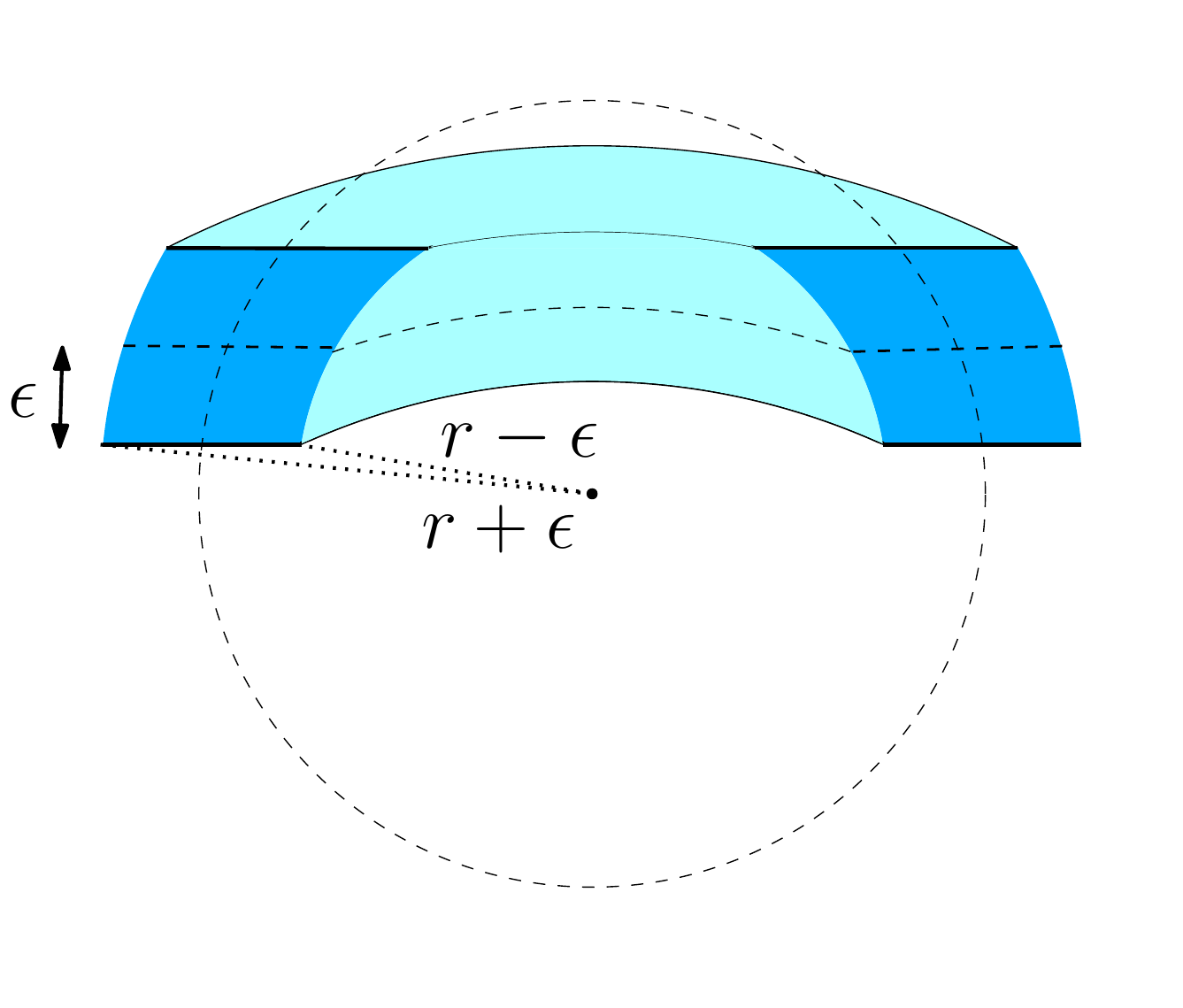}
\caption{\footnotesize{\textbf{Illustration of a slice of $\is(\epsilon)$:} in blue, a piece of $\is(\epsilon)(\sigma,\lambda,u)$. The slice $\is(\epsilon)(\sigma,\lambda,u)$ is obtained as the intersection of the $\epsilon$-thickening of the $(n-m+1)$-plane $P_u$ and the $\epsilon$-thickening of the $(n-m)$-sphere $S_u$ of radius $r$.}}\label{figure_slice_is_eps}
\end{figure}

From this expression it is readily seen that $g=L_v^*h$ if and only if $g(u,u) = h(v,v)$ and $v \in L(u) + L(\ker \lambda )^{\perp}$. So $v$ lies inside the $(n-1)$-dimensional sphere $S_u$ of radius $\|u\|_{g}$ and inside the affine $(n-m+1)$-plane
\begin{eqnarray*}
P_u := L(u) + L(\ker \lambda)^{\perp}.
\end{eqnarray*}
Thus $\is(\sigma,\lambda,u) = S_u \cap P_u$ is a $(n-m)$-dimensional sphere of $T_yW$ and its convex hull is a ball of the same dimension (see Figure~\ref{image_relation_isometrie}). Since we have assumed $n>m$, the space $\is(\sigma,\lambda,u)$ is arc-connected. Since $\is(\sigma,\lambda,u)$ is a $(n-m)$-dimensional sphere, $\iconv \, \is(\sigma,\lambda,u)$ is a $(n-m+1)$-ball of~$P_u$.\\

So a slice of the relation of $\epsilon$-isometric maps is a thickening of $\is (\sigma,\ell,u)$ (see Figure~\ref{figure_slice_is_eps}).


\subsection{Characterization of subsolutions of the relation of isometric maps}

Let $\proj_0$ be the orthogonal projection on $\ker \lambda$ in $T_xM$ and $\proj_P$ be the orthogonal projection on $P= L(\ker \lambda)$ in $T_yW$. We characterize subsolutions of $\is$ with respect to $(d\pi,u)$, for a submersion $\pi:U\subset M \rightarrow \R$ and a tangent vector field $u:U\rightarrow TM$ such that $d\pi(u) =1$, in the following proposition:

\begin{prop}\label{prop_subsolution_is}
Let $f_0:M\rightarrow W$ be a $C^1$-map and $P:=df_0(\ker d\pi)$ such that $\dim\; P(x)=m-1$ for all $x\in U$. If $f_0$ satisfies $g|_{\ker d\pi} = f_0^*h |_{\ker d\pi}$, then a section 
$$x\mapsto\s(x) = (x,f_0(x), L_x :=(df_0)_x + (v_x - (df_0)_x(u_x)) \otimes d\pi_x )$$
is a formal solution of $\is$ with respect to $(d\pi,u)$ if and only if, for every $x$, the vector $v_x$ can be written in the form $v_x =  \proj_{P(x)} L_x(u_x) + \tau_x$ where $\tau_x\in P(x)^{\perp}$ and $\|\tau_x\|_h = r(x) = \sqrt{\|u_x\|_{g}^2 - \| \proj_{0} u_x\|_{g}^2 }$. 
\end{prop}

\noi
\textbf{Proof.--} Recall that $v_x \in \is(\s(x), d\pi_x, u_x)$ if and only if $v_x\in  S_{u(x)} \cap P_{u(x)}$ i.e.
\begin{eqnarray*}
\|v_x\|_{h}^2 = \| u_x \|_{g}^2 \,\, \mbox{ and } \,\, \proj_{P(x)} v_x = \proj_{P(x)} L(u_x).
\end{eqnarray*}
Decomposing $v_x$ in $P(x)\oplus P(x)^{\perp}$, we have $v_x = \proj_{P(x)} L(u_x) + \widetilde{\tau_x}$, where $\widetilde{\tau_x}$ is a vector of $P(x)^{\perp}$ of norm $\|\widetilde{\tau_x}\|_h = r(x)$ by definition of $r$. Now we have to give an expression of the radius $r$ which only depends on $u$ and not to $\s$. By the Pythagorean theorem we have
\begin{eqnarray*}
r(x)^2 = \|\widetilde{\tau_x}\|_h^2 = \|v_x\|_h^2 - \|\proj_{P(x)} L(u_x)\|_h^2.
\end{eqnarray*}
As $\|v_x\|_h=\|u_x\|_g$, we then have $\|\widetilde{\tau_x}\|_h^2 = \|u_x\|_g^2 - \|\proj_{P(x)} L(u_x)\|_h^2$. The space $P=L(\ker \lambda)$ depends on $L$, so $\s$. Let $u_x = \proj_0 u_x +(u_x-\proj_0 u_x)$ with $\proj_0$ the orthogonal projection on $\ker \lambda$. Then
\begin{eqnarray*}
L(u_x) = L \Big( \proj_0 u_x +(u_x-\proj_0 u_x) \Big)=  L(\proj_0 u_x) + L(u_x - \proj_0 u_x).
\end{eqnarray*}
As $L$ is isometric we have, for any $a\in \ker \lambda$ and $b\in (\ker \lambda)^{\perp}$,
\begin{eqnarray*}
\langle a,b \rangle =0 \Leftrightarrow \langle L(a),L(b) \rangle =0.
\end{eqnarray*}
In particular, for $b= u - \proj_0 u$, that implies $L(u - \proj_0 u) \in L(\ker \lambda)^{\perp} = P(x)^{\perp}$. Thus $\proj_{P(x)} L(u_x) = L(\proj_0 u_x)$ and
\begin{eqnarray*}
 \| \proj_{P(x)} L(u_x) \|_h =  \| L(\proj_{0} u_x ) \|_h  =  \| \proj_{0}u_x \|_g
\end{eqnarray*}
the last equality comes from $L$ is isometric. So
\begin{eqnarray*}
r(x)^2 = \|\widetilde{\tau_x}\|_h^2 = \|v_x\|_h^2 -  \| \proj_{0}u_x \|_g^2.
\end{eqnarray*}
\hfill $\square$

\subsection{Proof of Theorem~\ref{th_is_Kuiper}}

We begin with a preparatory lemma, then describe $\iconv\, \is(\epsilon)(\sigma,\lambda,u) \cap P_u(w)$ and define a $c$-shaped  loop family for the relation $\is(\epsilon)$. We finally construct $\bgamma$ and prove that it is surrounding.\\

Let $\sigma=(x,y,L)\in \is(\epsilon)$. Let $\lambda\in T_x^*M$, $u\in T_xM$ such that $\lambda(u) =1$, and let $w \in \iconv\, \is(\epsilon)(\sigma,\lambda,u)$. Note that as $\is(\epsilon)$ is a thickening of $\is$ and by definition of $\sigma$ and $w$, the distance (for the metric $h$) between $w$ and $P_u$ is less than $2\epsilon$, but $w$ does not belong necessarily to $P_u$. We denote by $P_u(w)$ the affine $(n-m+1)$-plane that contains $w$ and which is a translation of $P_u$:
$$P_u(w) := \{ v \in T_yW \;|\; \proj_{P} w = \proj_{P} v \}$$
where $P$ denotes $L(\ker \lambda)$. Thanks to the following lemma, we can assume that $w$ belongs to $P_u$:
\begin{lemme}\label{lemme_is_homotopie}
Let $(\sigma,w) \in \iconv(\is(\epsilon), \lambda, u)$ with $\sigma=(x,y,L)$. There exists a homotopy $\sigma_t = (x,y,L_t)$ such that $\sigma_0 = \sigma$, $\sigma_t\in\is(\epsilon)(\sigma,\lambda,u)$ for all $t\in [0,1],$ and $\proj_{P} L_1(u) = \proj_{P} w$.
\end{lemme}

\noi
\textbf{Proof.--} We set $v_0 = L_0(u) = L(u)$. We can assume that $\|v_0\|_h \geq \|w\|_h$. Indeed, if $\|v_0\|_{h} < \|w\|_{h}$ we perform a first homotopy. Let $\widetilde{L_t} = L +(\widetilde{v_t} - v_0)\otimes \lambda$ where 
\begin{eqnarray*}
\widetilde{v_t} := \proj_{P} v_0 + \left( (1-t) + t\frac{\sqrt{\|w\|_h^2-\|\proj_{P} v_0\|_h^2}}{ \|v_0 - \proj_{P} v_0\|_h}\right)( v_0-\proj_{P} v_0).
\end{eqnarray*}
This homotopy joins $v_0$ to $\widetilde{L_1}(u)=\widetilde{v_1}$ where $\|\widetilde{v_1} \|_h = \|w\|_{h}$. Let $V_0 = v_0$ if $\|v_0\|_{h} \geq \|w\|_{h}$, and $V_0 =\widetilde{v_1}$ if $\|v_0\|_{h} < \|w\|_{h}$. In both cases, we consider the homotopy $L_t = L +(v_t - V_0)\otimes \lambda$ with:
\begin{eqnarray*}
v_t := t\, \proj_{P} w + (1-t)\, \proj_{P} V_0 + \varphi(t) (V_0 - \proj_{P} V_0)
\end{eqnarray*}
and
\begin{eqnarray*}
\varphi(t) = \sqrt{ \frac{\| V_0 \|_h^2-\|t\, \proj_{P} w +(1-t)\, \proj_{P} V_0 \|_h^2}{\|V_0 - \proj_{P} V_0 \|_h^2} }.
\end{eqnarray*}
Since $\|V_0\|_{h} \geq \|w\|_{h}$ the numerator is positive and $\varphi$ is well defined. By definition of $\varphi$, for every $t$, we have $\|v_t\|_{h} = \|V_0\|_{h}$. This property ensures that $\sigma_t=(x,y,L_t)\in\is(\epsilon)(\sigma,\lambda,u)$ for all $t\in [0,1].$ By the expression of $v_t$, we have $\proj_{P} v_1 =\proj_{P} w$. \hfill $\square$\\

This lemma and Point $(3)$ of Definition~\ref{def_surrounding} imply that it is enough to construct the loop family $\bgamma$ for every couple $(\sigma,w)$ such that $\proj_{P} L(u) =\proj_{P} w$. We assume in the sequel that this last condition is fulfilled together with the fact that the codimention is one.\\

\noi
\textbf{Description of $\iconv\, \is(\epsilon) (\sigma,\lambda,u) \cap P_u(w)$.--} By assumption $n=m+1$ therefore the space $P_u(w)$ is a 2-plane. We denote by $D(\rho)$ the open disk of $P_u(w)$ with radius $\rho$ and center $\mbox{proj}_P(L(u))$ and by $A(\rho_{min},\rho_{max})$ the open annulus $D(\rho_{max})\setminus \overline{D(\rho_{min})}.$ The intersection of the thickened relation $\is(\epsilon)(\sigma,\lambda,u)$ with $P_u(w)$ is either an annulus or a disk depending on the value of $\epsilon$. Precisely, let 
\begin{eqnarray*}
\begin{array}{lll}
r_{min}^2(\epsilon)& := & \displaystyle \min\left((\|u\|_{g} - \epsilon)^2-\|\proj_P w \|_h^2,0\right)\\
r_{max}^2(\epsilon)& := & \displaystyle (\|u\|_{g} + \epsilon)^2-\|\proj_P w\|_h^2.
\end{array}
\end{eqnarray*}
because the sphere $S_u$ of Paragraph~\ref{subsect_geom_desc} is of radius $\|u\|_{g}$. A computation shows that $\is(\epsilon)(\sigma,\lambda,u)\cap P_u(w)$ is the annulus $A(r_{min}(\epsilon),r_{max}(\epsilon))$ if $r_{min}(\epsilon)>0$ and the disk $D(r_{max}(\epsilon))$ if $r_{min}(\epsilon)=0$. In any case,
\begin{eqnarray*}
\iconv \,\is(\epsilon)(\sigma,\lambda,u) \cap P_u(w) = D(r_{max}(\epsilon)).
\end{eqnarray*}
In particular, we have $w\in D(r_{max}(\epsilon))$ and $L(u)\in A(r_{min}(\epsilon),r_{max}(\epsilon))$. We want to build a $c$-shape loop family inside $A(r_{min}(\epsilon),r_{max}(\epsilon))$, for that we define a disk which will support $\bgamma$ and such that a neighborhood of this disk will be in $A(r_{min}(\epsilon),r_{max}(\epsilon))$ too. Let $D(\widetilde{r})$ a disk where
\begin{eqnarray*}
\widetilde{r} = \max \left( \sqrt{ \|L(u)\|_h^2 - \|\proj_{P} L(u)\|_h^2 } , \sqrt{\|w\|_h^2-\|\proj_{P} w\|_h^2}+\frac{1}{3}d_1(w)\right)
\end{eqnarray*}
where
\begin{eqnarray*}
d_1(w) := dist(w,\partial (\iconv\, \is(\epsilon)(\sigma,\lambda,u)) )
\end{eqnarray*}
is the distance between $w$ and the boundary of the convex hull of $\is(\epsilon)(\sigma,\lambda,u)$. Moreover we have $w\in D(\widetilde{r})$ and $\partial D(\widetilde{r})\subset A(r_{min}(\epsilon),r_{max}(\epsilon)).$\\

\noi
\textbf{Parametrization of $D(\widetilde{r})$.--} Let $\nu$ be the unique unit normal vector of $L(T_xM)$ induced by the orientation of $M$ and $W$. We see $P_u(w)$ as the complex plane $\C$ by identifying the base $\left(\nu ,(L(u)-\proj_P L(u))/ \|L(u)-\proj_P L(u)\|_h\right)$ with $(1,i)$ and we define a parametrization of $\overline{D(\widetilde{r})}$ by 
\begin{eqnarray*}
\begin{array}{lcll}
b : &  [0,\pi]\times [0,1] & \longrightarrow & \overline{D(\widetilde{r})}\\
& (\theta,\beta) &  \longmapsto & \proj_P L(u)+ \beta \widetilde{r} e^{i\theta} + (1-\beta) \widetilde{r} e^{-i\theta}.
\end{array}
\end{eqnarray*}
This parametrization is 1-to-1 except over points of the form $(0,\beta)$ and $(\pi,\beta)$. It maps the boundary of the square  $[0,\pi]\times [0,1]$ onto the circle $\partial D(\widetilde{r})$.\\

\begin{figure}[!ht]
\centering
\includegraphics[scale=0.64]{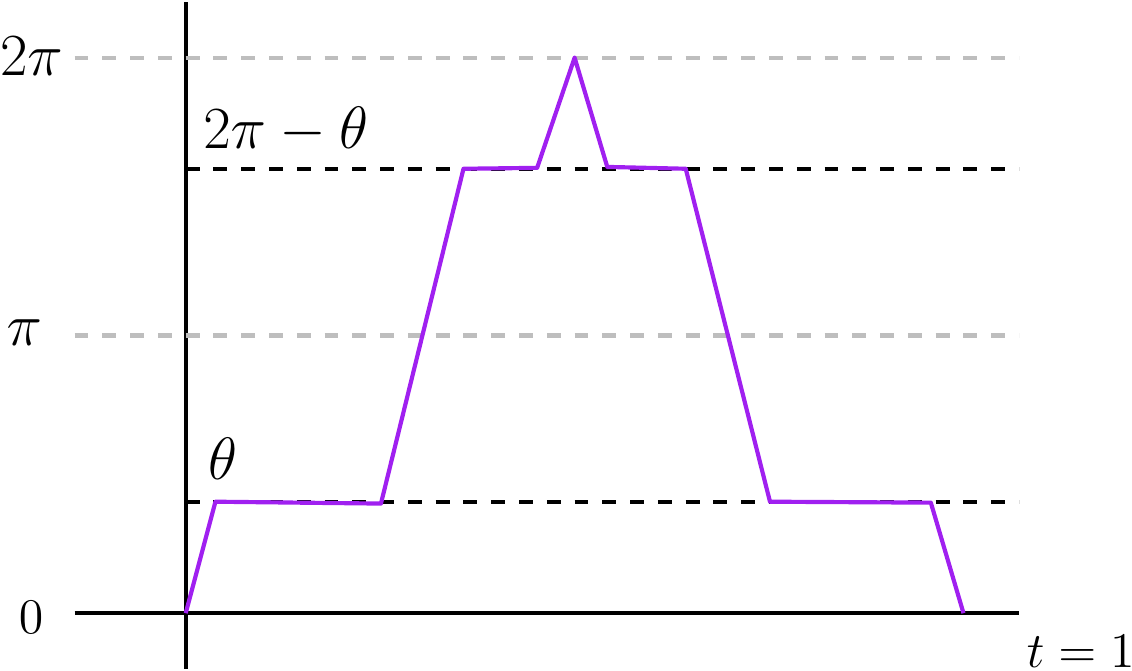}
\includegraphics[scale=0.66]{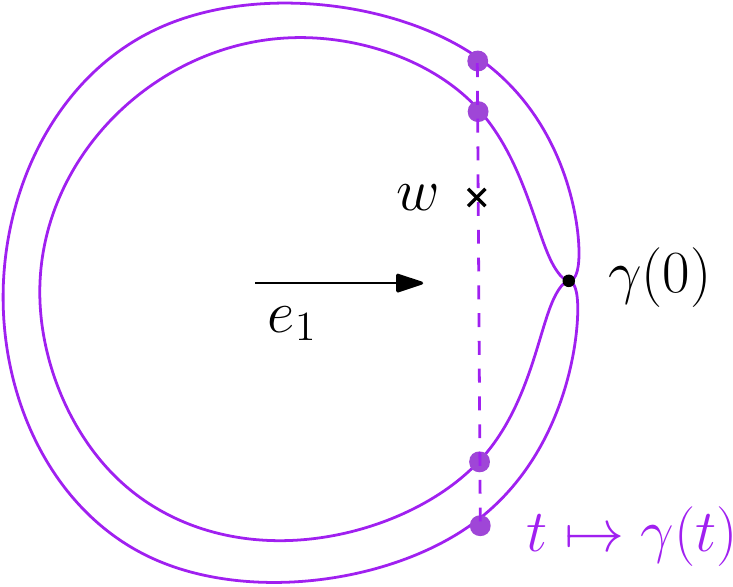}
\caption{\footnotesize{\textbf{Proof of Theorem~\ref{th_is_Kuiper}}: Left: the graph of the function $g_{\theta, \beta}$, Right: the image of the loop $\gamma$ in the affine plane $P_u(w)$, the two circles visualise the round-trip of the loop.}} \label{image_shape_isom}
\end{figure}

\noi
\textbf{The shape.--} We first define the parameter space $A$ to be
\begin{eqnarray*}
A := \{ (\eta,\theta,\beta)\in\; ]0, \textstyle \frac{1}{2}[\times [0,\pi]\times [0,1]\; |\; \eta\leq \beta\leq 1-\eta\}.
\end{eqnarray*}
and then the shape  $c:A\times \R/\Z\rightarrow \C\times\R$ by
\begin{eqnarray*}
c(\eta,\theta, \beta,t):=(\exp (i g_{\theta,\beta}(t))+\eta \cos\theta, 1).
\end{eqnarray*}
The image of $c(\eta,\theta,\beta, \cdot)$ is a whole circle of center $(\eta\cos\theta,1)$ and radius 1. Let $\beta' = \beta - \frac{\eta}{2}$, the angular function $g_{\theta,\beta}$ is the piecewise linear map given by
\begin{itemize}
\item[(i)] $g_{\theta,\beta}(0)=0$ and $g_{\theta,\beta}(\frac{1}{2})=2\pi$
\item[(ii)] $g_{\theta,\beta}(t) = \theta 
\mbox{ on } \displaystyle \left[ \frac{\eta\theta}{4\pi}, \frac{\beta'}{2} + \frac{\eta \theta}{4\pi} \right]$
\item[(iii)] $g_{\theta,\beta}(t) =  2\pi - \theta 
\mbox{ on } \displaystyle \left[ \frac{\beta'}{2} + \frac{\eta(2\pi - \theta)}{4\pi}, \frac{1}{2} -\frac{\eta \theta}{4\pi} \right]$
\end{itemize}
on $[0,\frac{1}{2}]$ and such that $g_{\theta,\beta}(t) = g_{\theta,\beta}(1-t)$ for all $t\in [0,\frac{1}{2}]$ (see its graph on Figure~\ref{image_shape_isom}). A computation shows that 
\begin{eqnarray*}
\overline{ c(\eta,\theta,\beta) } 
&=&  (\beta e^{i\theta} + (1-\beta)e^{-i\theta},1).
\end{eqnarray*}

\noi
\textbf{The loop family.--} Since $b$ induces a bijection between $]0,\pi[\times ]0,1[$ and $D(\widetilde{r})$, there exists a unique couple $(\theta,\beta)\in\; ]0,\pi[\times ]0,1[$ such that $b(\theta,\beta)=w$. We define two functions $c_1$ and $c_2$ by the equality
$$c(\eta,\theta,\beta,\cdot)=(c_1(\cdot)+ic_2(\cdot),1)$$
($\eta$ will be chosen later). We put 
\begin{eqnarray*}
\textbf{e}_1 := \widetilde{r} \frac{\nu}{\|\nu\|_h}, \;\;\;  \textbf{e}_2 := \widetilde{r} \frac{ L(u) - \proj_P  L(u) }{\|L(u) - \proj_P L(u)\|_h}\;\;\;\mbox{and}\;\;\; \textbf{e}_3 := \proj_P  L(u)
\end{eqnarray*}
and we define the loop family $\bgamma$ by
\begin{eqnarray*}
\bgamma(\sigma,w)(t) := c_1(t) \textbf{e}_1 + c_2(t) \textbf{e}_2 + \textbf{e}_3.
\end{eqnarray*}

\begin{figure}[!ht]
\centering
\includegraphics[scale=0.7]{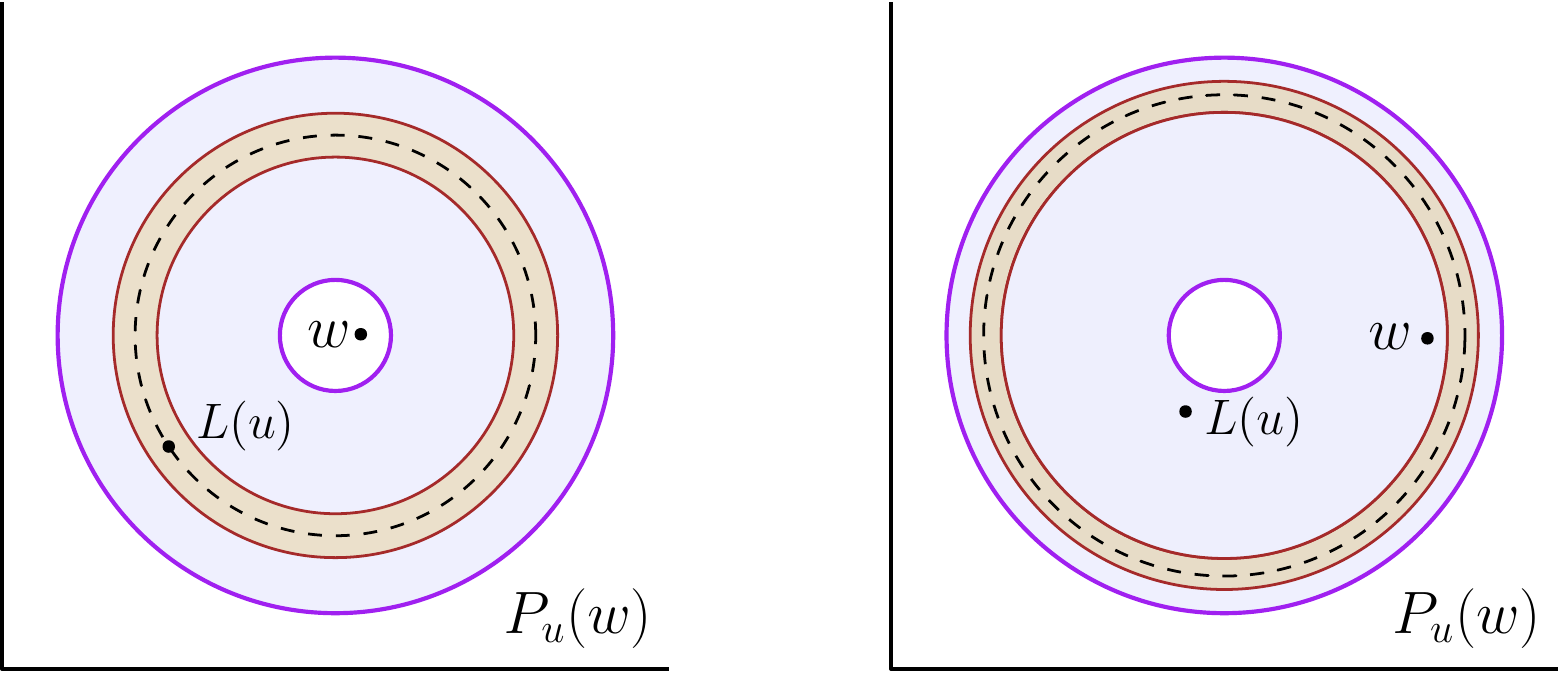}
\caption{\footnotesize{\textbf{Proof of Theorem~\ref{th_is_Kuiper}}: In purple, the trace of the relation on $P_u(w)$, in black (dashed line) the boundary of disk $D(\widetilde{r})$ and in brown, the annulus $A(\tilde{r}(1-\eta),\tilde{r}(1+\eta))$ depending on whether $\sqrt{\|L(u)\|_h^2-\|\proj_P L(u)\|_h^2} > \sqrt{\|w\|_h^2-\|\proj_P w\|_h^2}+\frac{1}{3}d_1(w)$ (left) or not (right), see the definition of $\widetilde{r}$.}}
\end{figure}

\noi
The image of the loop $\bgamma(\sigma,w)$ is the translated circle $\partial D(\tilde{r})+\tilde{r}\eta\cos\theta\textbf{e}_1$ which lies inside the annulus $A(\tilde{r}(1-\eta),\tilde{r}(1+\eta))$ of $P_u(w)$. Consequently, to ensure that the image of $\bgamma(\sigma,w)$ is in the relation, it is enough to choose $\eta$ such that $A(\tilde{r}(1-\eta),\tilde{r}(1+\eta))\subset A(r_{min}(\epsilon),r_{max}(\epsilon))=\is(\epsilon)(\sigma,\lambda,u)\cap P_u(w)$. It is readily checked that the choice 
$$\eta := \frac{1}{3}\min(d_1(w),d_2(L(u)))$$
where $d_2(L(u)) = dist(L(u),\partial A(r_{min}(\epsilon),r_{max}(\epsilon)))$ is convenient. It is also straightforward to see that this loop family satisfies the Average Constraint: $\overline{\bgamma}(\sigma,w) = w.$ The base point of the loop is $\bgamma(\sigma,w)(0) = (1+\eta\cos\theta)\textbf{e}_1 +\textbf{e}_3$. The homotopy $H(s):=(\cos s\;\textbf{e}_1+\sin s\;\textbf{e}_2) +\eta\cos\theta\textbf{e}_1+\textbf{e}_3$ with $s\in\; [0,\frac{\pi}{2}]$ connects $\bgamma(\sigma,w)(0)$ with $\eta\cos\theta\textbf{e}_1 +\textbf{e}_2+ \textbf{e}_3$
A linear homotopy joins this last point to $L(u) = (\|L(u) - \proj_P L(u) \|_h )\textbf{e}_2 / \widetilde{r} + \textbf{e}_3.$ Consequently, the loop family $\bgamma$ is $c$-shaped and surrounding (see Definition~\ref{def_surrounding}). This proves that $\is(\epsilon)$ is a Kuiper relation.

\section{An application: desingularization of a cone to a surface $\epsilon$-isometric to a flat cylinder}

Proposition~\ref{prop_noInt} together with the Kuiper property  of the relation of $\epsilon$-isometric maps are the reason of the absence of integrals in the formula proposed in \cite{Kuiper55,ContiEtAl} to solve $\is(\epsilon)$. The approach developed here also allows to apply the $h$-principle in its full generality for $\is(\epsilon)$. Indeed, in the above cited references, the formulas only make sense when the base map $f_0$ is an immersion but in the framework of the $h$-principle this hypothesis is not required: provided that $\s$ is a subsolution, any base map $f_0$, singular or not, is convenient.\\

Here, we illustrate this point with a basic example. We consider a singular map sending a flat cylinder onto a cone and we use the Kuiper property of $\is(\epsilon)$ to build an $\epsilon$-isometric map arbitrarily closed (in the $C^0$ sense) to the initial singular map.

\subsection{Formal solution}
We identify the flat cylinder of height $\frac{1}{20}$ and radius $\frac{1}{2\pi}$ with the space $Cyl=\R/\Z\times [-10^{-1},10^{-1}]$ endowed with the Euclidean metric. We define our formal solution to be 
$$\s_0 :(x,y) \longmapsto ((x,y), f_0(x,y), v_1(x,y), \partial_2 f_0(x,y))$$
where $f_0$ is a parametrization of a cone:
$$f_0(x,y) = \frac{1}{\sqrt{2}} \Big( y\cos(2\pi x ), y \sin(2\pi x), y \Big)$$
and $v_1$ is such that $\partial_1 f_0(x,y) = \sqrt{2}\pi  y v_1(x,y)$. Precisely:
$$v_1(x,y) = \Big(- \sin(2\pi x), \cos(2\pi x) , 0\Big). $$
Observe that for every $(x,y)$ we have
\begin{eqnarray*}
\left\{
\begin{array}{l}
\|v_1(x,y)\| =1\\
\|\partial_2 f_0(x,y)\| =1 \\
\langle v_1(x,y), \partial_2 f_0(x,y) \rangle =0
\end{array}
\right.
\end{eqnarray*}
so the section $\s_0$ is a formal solution of the relation of $\epsilon$-isometric maps for every $\epsilon>0$.

\begin{figure}[!ht]
\centering
\includegraphics[width=0.5\textwidth]{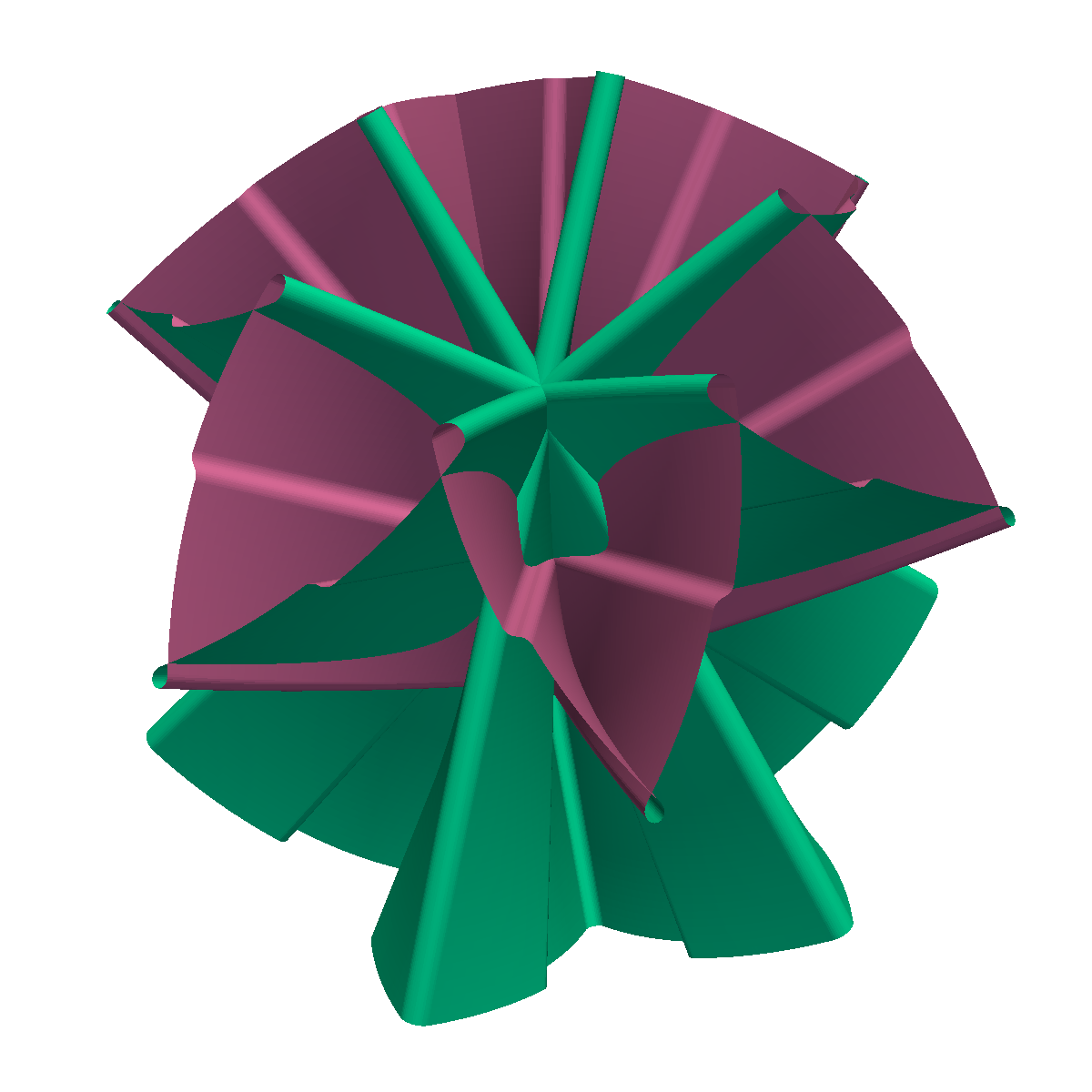}\includegraphics[width=0.5\textwidth]{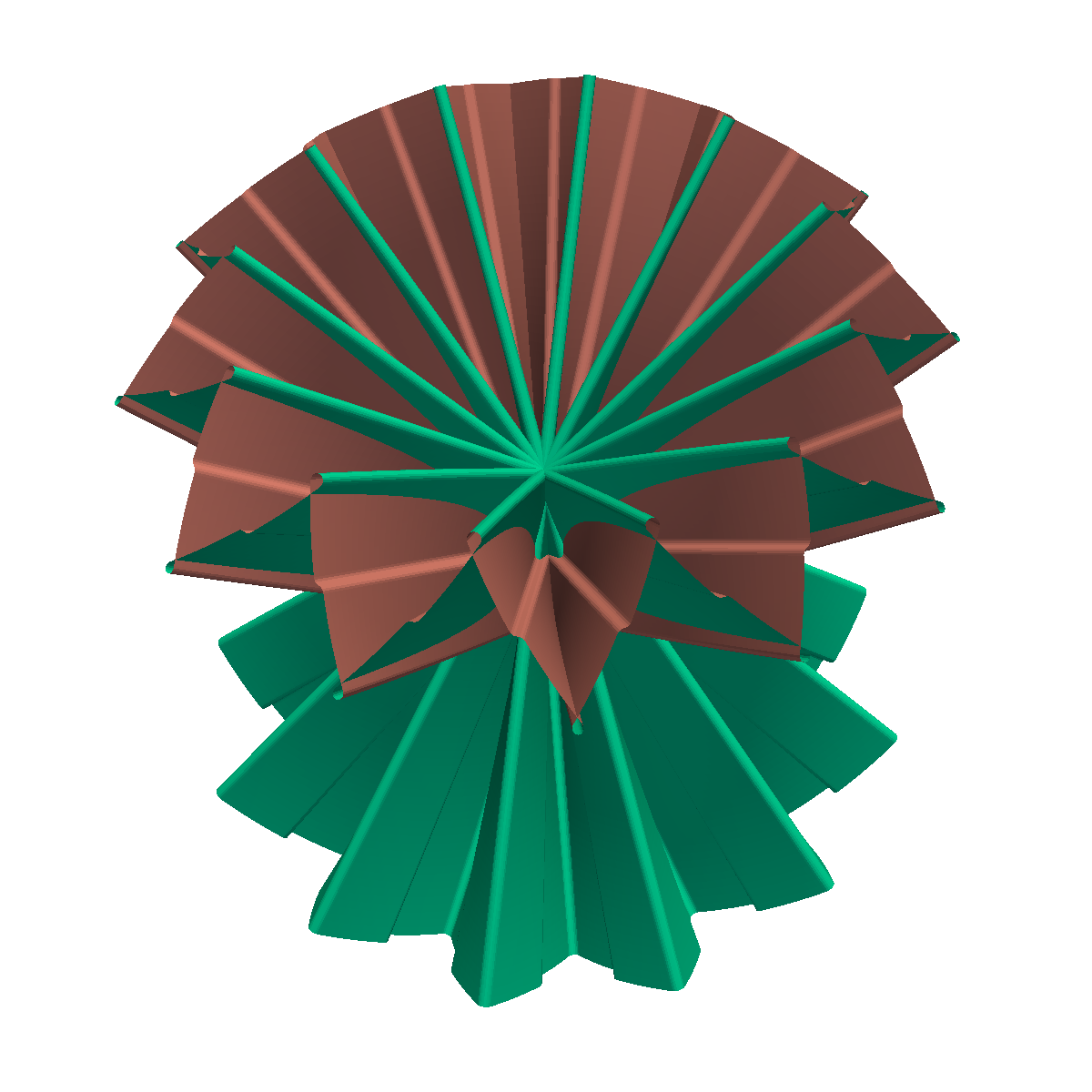}\\
\includegraphics[width=0.5\textwidth]{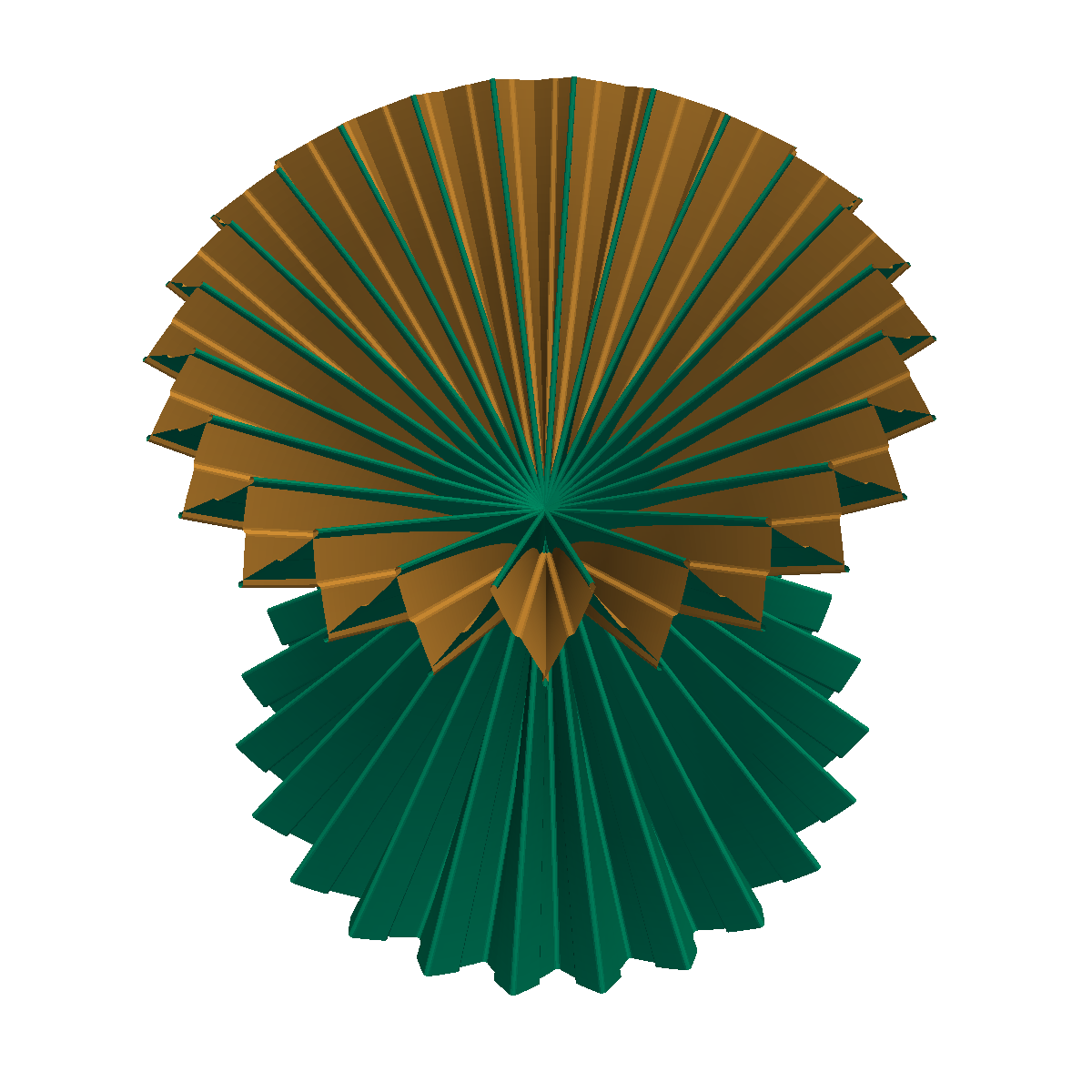}\includegraphics[width=0.5\textwidth]{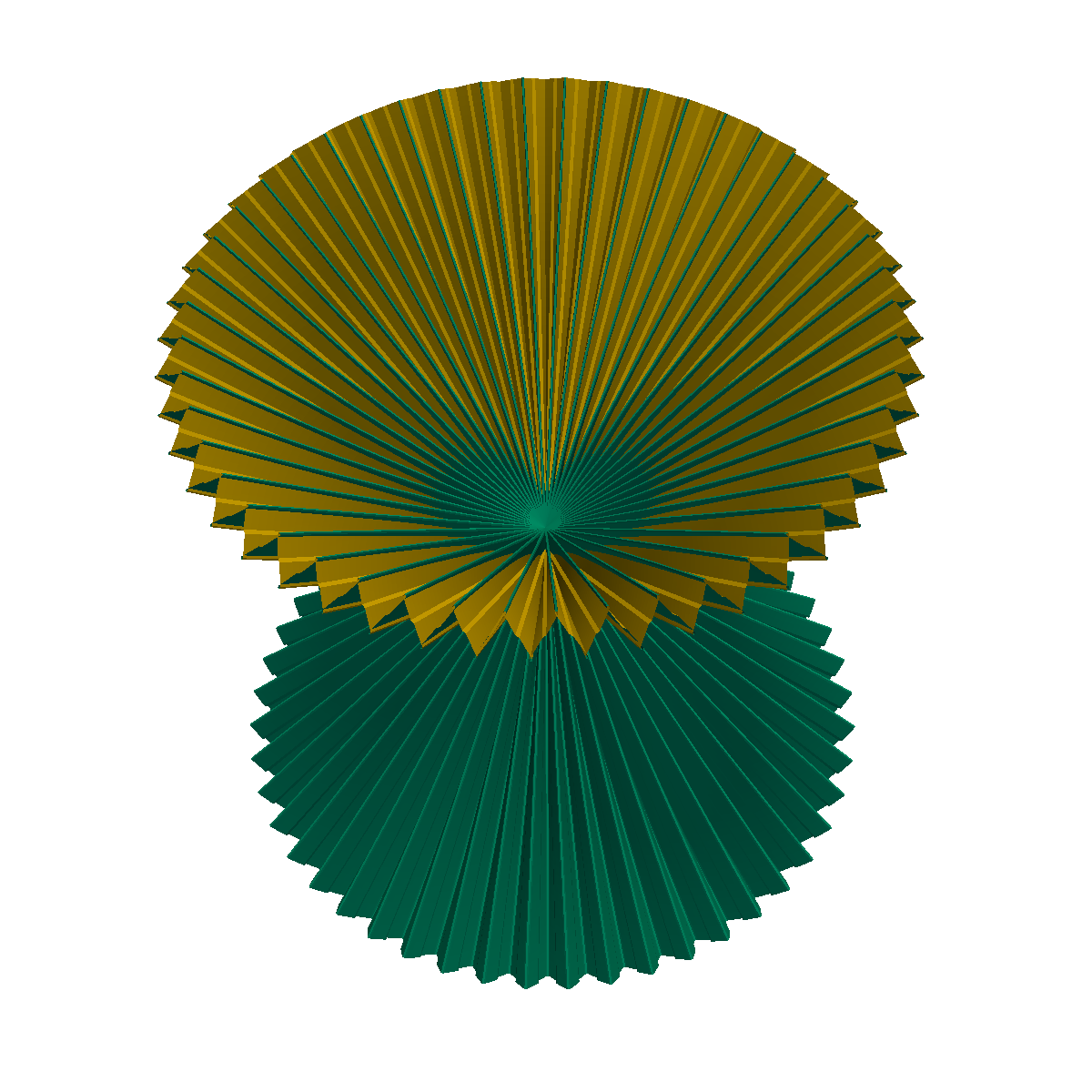}
\caption{\footnotesize{\textbf{Corrugation Process applied from a cone}}: Several images of $f_1(Cyl)$ with $\eta=0.2$ and, from left to right and up to down, $N=6,12,24,48$. Observe the $C^0$-density property (see Proposition~\ref{propCP} $(P1)$) : the larger $N$, the closer the surface to the cone.}
\end{figure}

\subsection{Subsolution}
The section $\s_0$ fails to be holonomic only in its $v_1$-component. To obtain a holonomic section, we thus intend to apply a Corrugation Process in the direction $\partial_1$. To do so, we need to check that $\s_0$ is a subsolution with respect to $\partial_1$. As $v_1$ and $\partial_2 f_0$ are orthogonal, the slice  $\is(\s_0,\partial_1)$ lies inside the plane spanned by $v_1$ and the normal vector
\begin{eqnarray*}
n(x,y) = v_1(x,y) \wedge \partial_2 f_0(x,y) = \frac{1}{\sqrt{2}} \Big( \cos(2\pi x), \sin(2\pi x), -1\Big)
\end{eqnarray*}
(see the proof of Theorem~\ref{th_is_Kuiper}). This slice is a circle of radius 1. The section $\s_0$ is a subsolution if and only if the derivative $\partial_1 f_0$ lies in the convex hull of $\is(\s_0,\partial_1)$. This condition is equivalent to
$|y| < (\sqrt{2}\pi)^{-1}.$
Since $y\in [-10^{-1},10^{-1}]$, this last inequality is fulfilled. This shows that $\s_0$ is subsolution with respect to $\partial_1$ of $\is(\s_0,\partial_1)$, and thus of $\is(\epsilon)(\s_0,\partial_1)$ for every $\epsilon>0.$

\begin{figure}[!ht]
\centering
\includegraphics[width=1\textwidth]{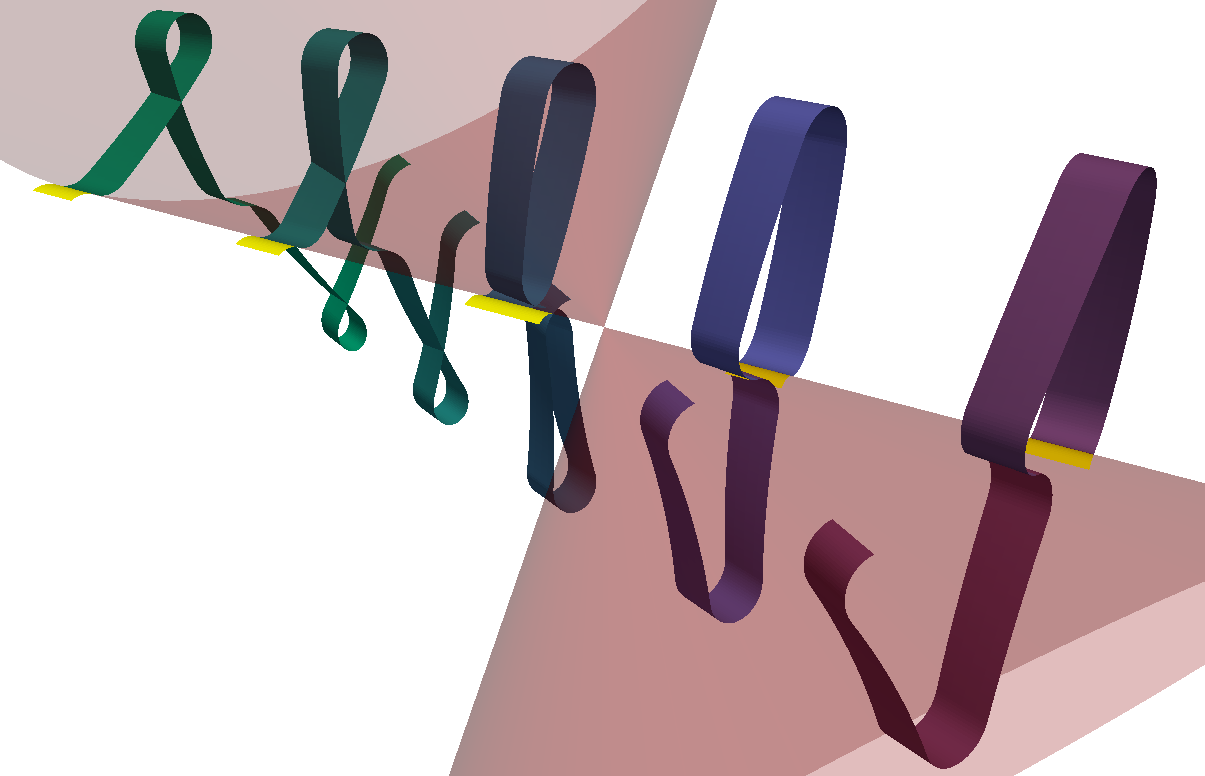}
\caption{\footnotesize{\textbf{The change of the shape of a corrugation when passing the conical singularity.}}: Here $N=6$ and $\eta=0.4$. }
\end{figure}

\begin{figure}[!ht]
\centering
\includegraphics[width=0.5\textwidth]{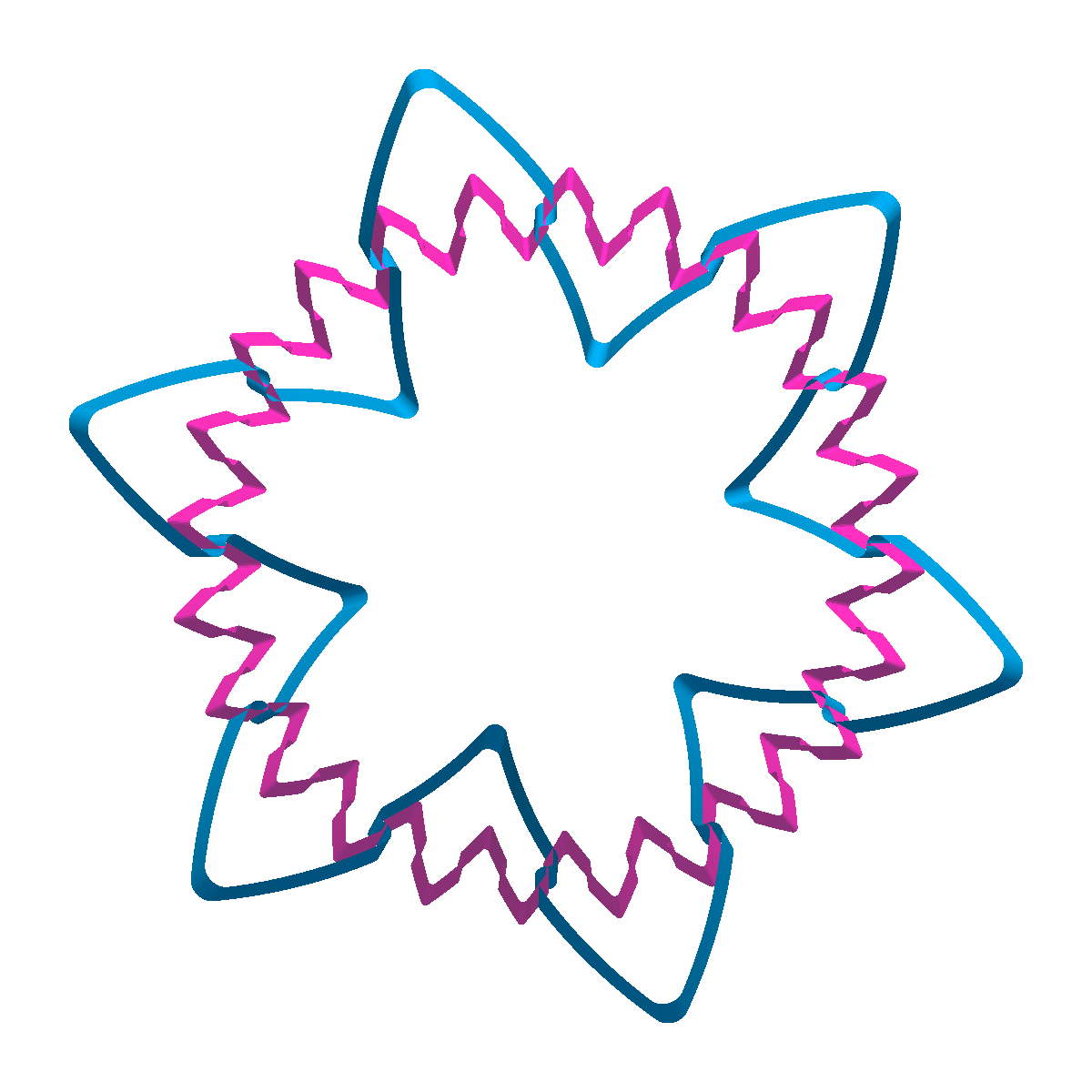}\includegraphics[width=0.5\textwidth]{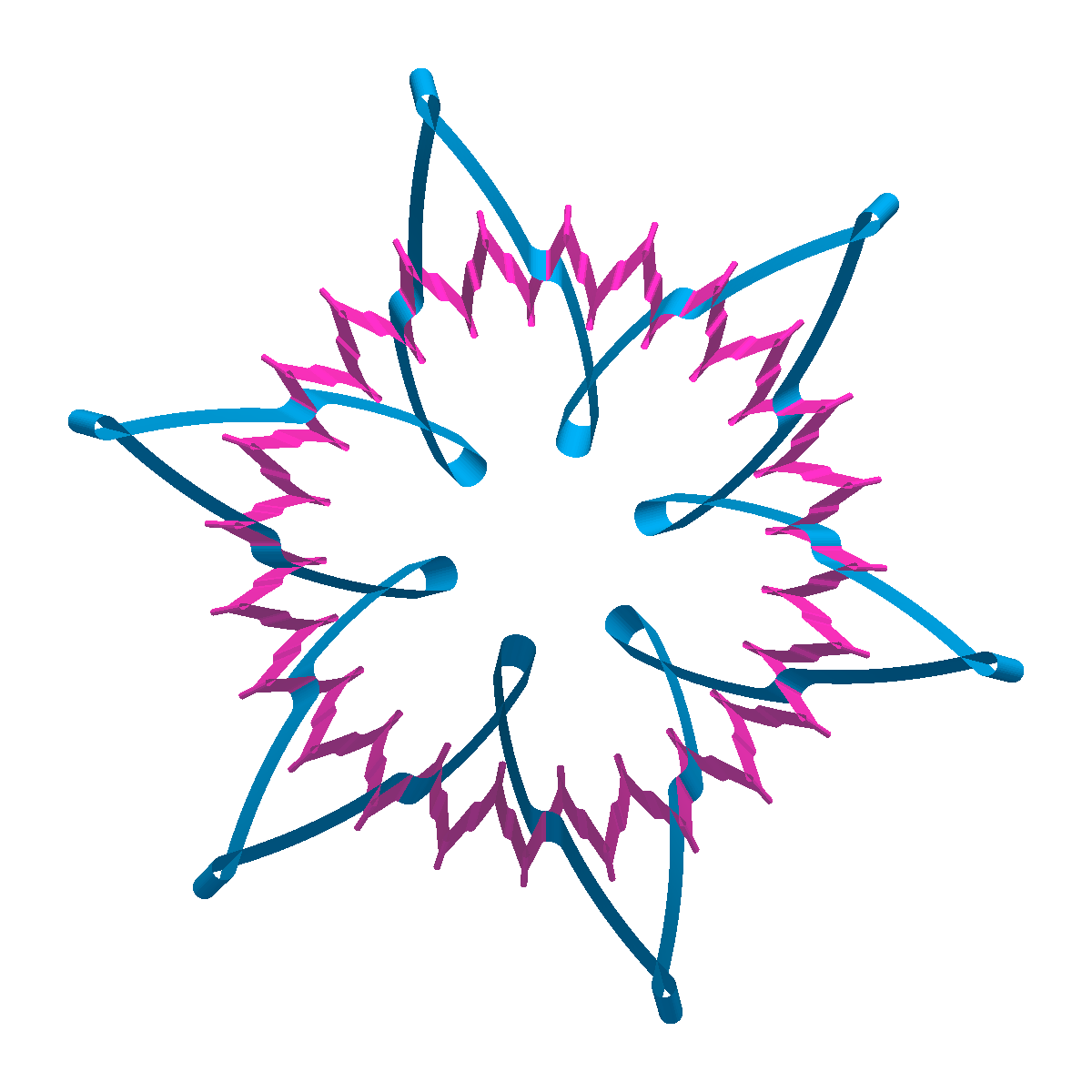}
\caption{\footnotesize{\textbf{Lengths and Corrugations}}: In blue, two slices of the image $f_1(Cyl)$ for $N=6$, in pink, two slices for $N=24$. On the right the slices are above the horizontal plane passing though the vertex of the cone. They are below this plane on the left. In all cases $\eta=0.2$.}
\end{figure}

\subsection{Corrugation Process}
We consider the shape $c:Cyl\times\R/\Z\rightarrow\C\times\R$ defined in subsection~\ref{th_is_Kuiper} where $\C$ is identified with the plane spanned by $(v_1,n)$:
\begin{eqnarray*}
c(\eta,\theta, \beta,t):=(\exp (i g_{\theta,\beta}(t))+\eta \cos\theta, 1).
\end{eqnarray*}
In that expression, $\beta$ and $\theta$ are defined by the relation
$$\partial_1 f_0=\beta e^{i\theta}+(1-\beta)e^{-i\theta}.$$

\begin{figure}[!ht]
\centering
\includegraphics[width=0.5\textwidth]{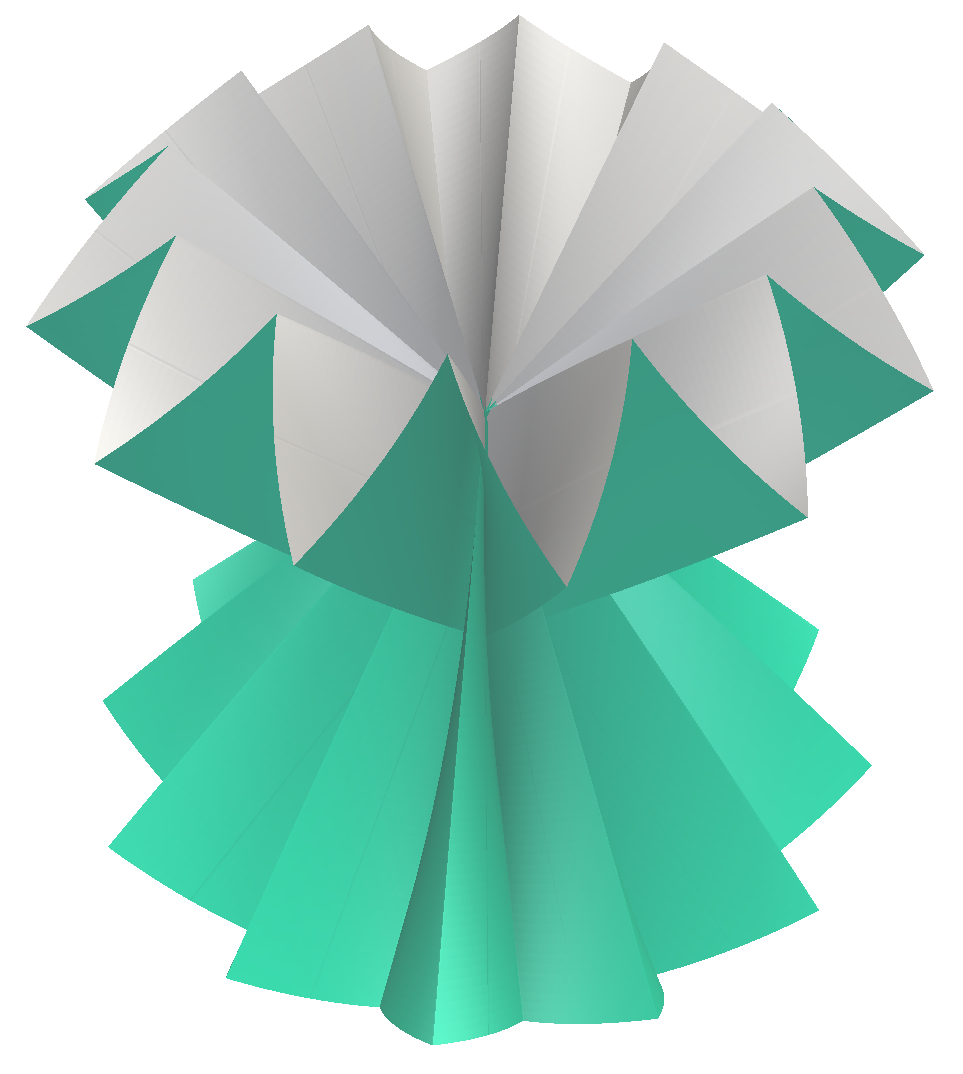}\includegraphics[width=0.5\textwidth]{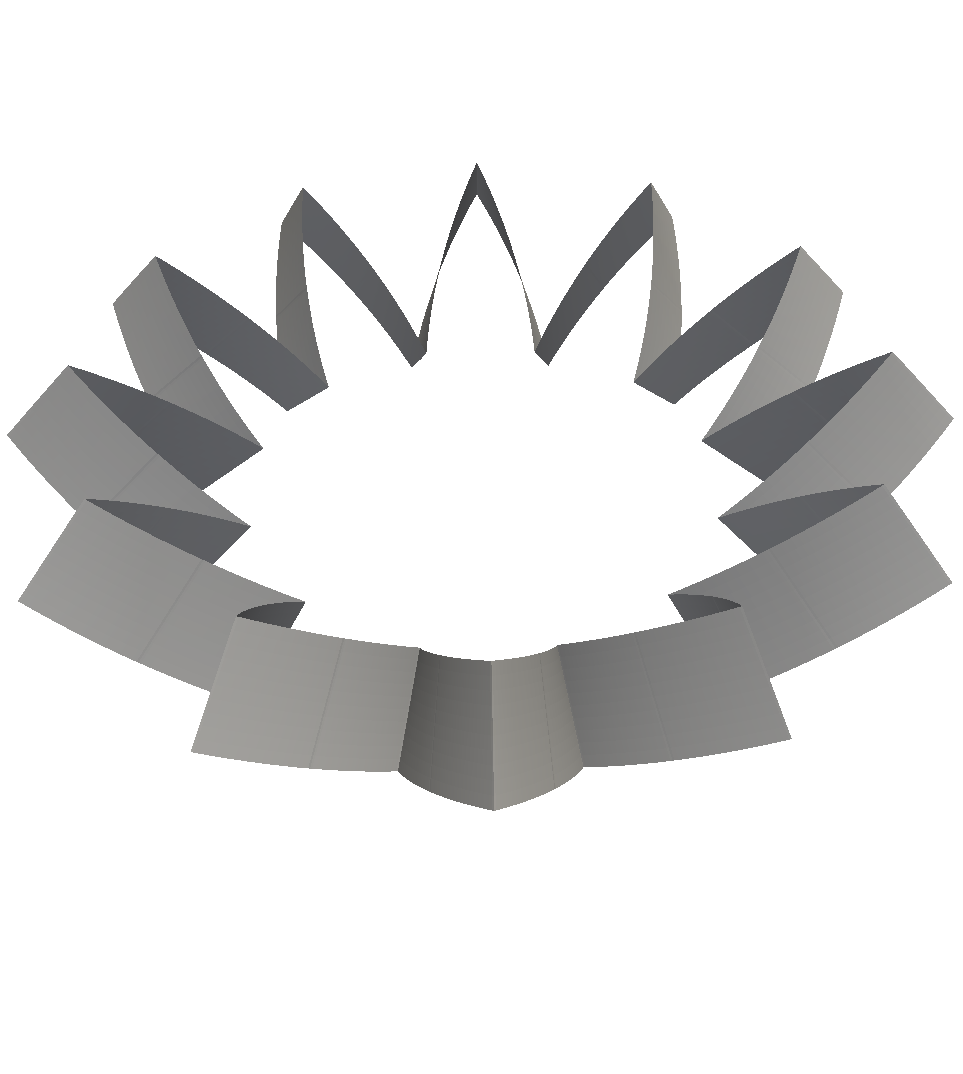}\\
\includegraphics[width=\textwidth]{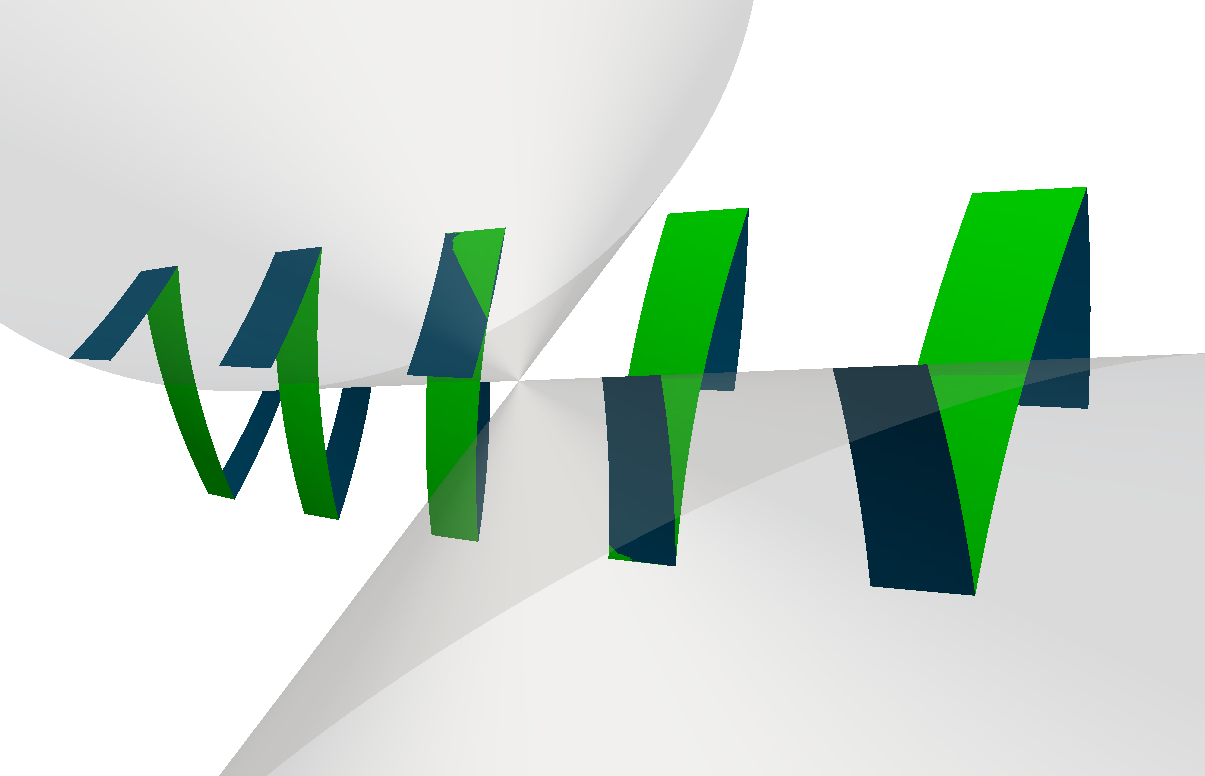}
\caption{\footnotesize{\textbf{PL behavior}}: Here $\eta=10^{-3}$ and $N=12$. Despite appearances, the map $f_1$ is still a $C^1$-immersion. See the close up of Figure~\ref{zoom_pl_behavior}.}\label{pl_behavior}
\end{figure}

Since $v_1$ is collinear to $\partial_1 f_0$ the coefficient $\beta$ is constant equal to $1/2$ and $\theta = \arccos (\langle \partial_1 f_0, v_1\rangle)= \arccos (\sqrt{2}\pi y) $. For short we denote $g$ for $g_{\theta,\beta}.$ The loop family $\gamma$ is thus given by
$$ \gamma(x,y,t) = \Big( \cos(  g(x,y,t)) + \eta \|\partial_1f_0(x,y)\|\Big) v_1(x,y) + \sin(x,y,t) n(x,y).$$
Observe that $\gamma(x,y,t)\in\is(\eta)(\s_0,\partial_1)$.
The Corrugation Process generates a map $f_1(x,y)= f_0(x,y)+ \frac{1}{N} \Gamma(x,y,Nx)$ with
\begin{eqnarray*}
\Gamma(x,y,t) := \int_{s=0}^t \gamma(x,y,s) - \overline{\gamma(x,y)}ds.
\end{eqnarray*}
Recall that, from Point $(P_3)$ of Proposition~\ref{propCP},
we have 
$$\partial_{1} f_1(x,y) = \gamma(x,y,Nx) + O(1/N).$$
Let $\epsilon>0$ be given. To insure $\partial_1f_1\in\is(\epsilon)(\s_0,\partial_1)$ we have to choose $\eta<\epsilon$ and $N$ large enough. 

\begin{figure}[!ht]
\centering
\includegraphics[width=\textwidth]{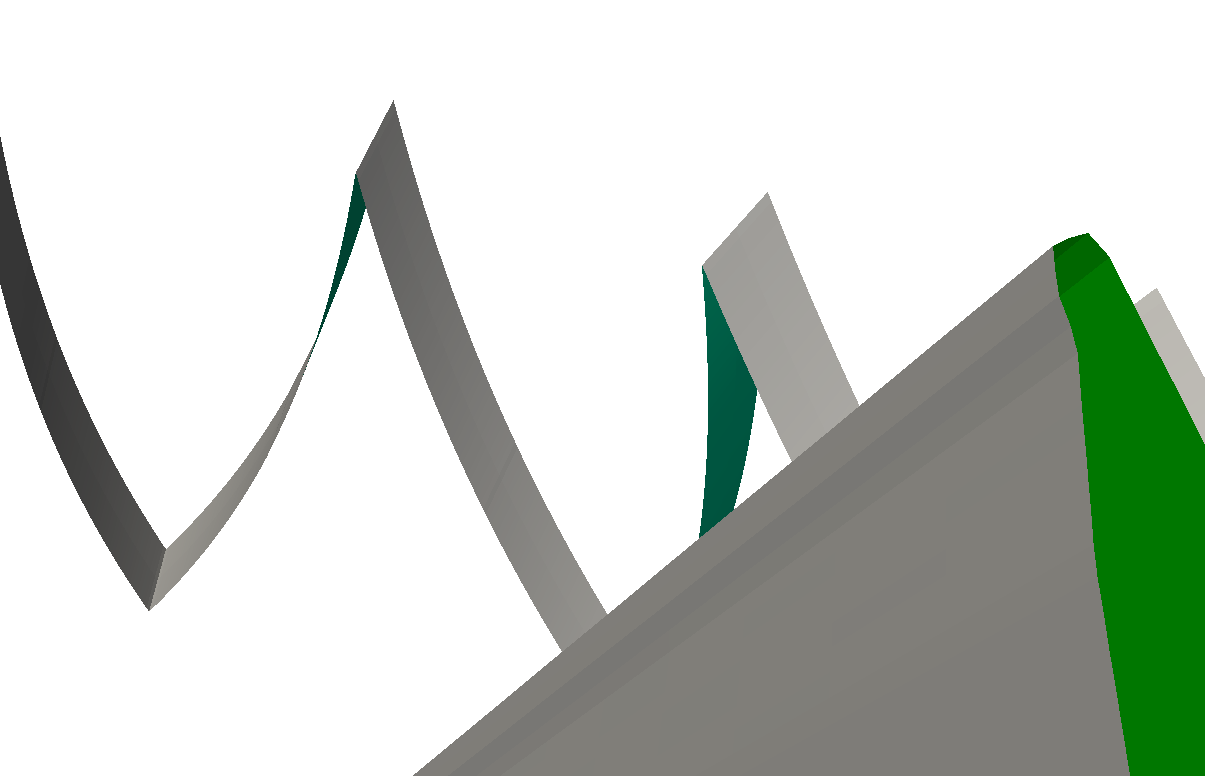}
\caption{\footnotesize{\textbf{Zoom on the peak of a corrugation}}: The peak is not a folding. For $N$ large enough, the corrugations are immersed. A close-up shows the roundness of the peak (in the foreground). The angles that appear are artefact due to the discretisation step. }\label{zoom_pl_behavior}
\end{figure}

\subsection{Numerical implementation}

We use the analytical expression of Proposition~\ref{prop_noInt} together with the above expression of $\gamma$ to implement the Corrugation Process. The images reveal corrugations whose shape varies from a small loop to the one of a roof. A closer look to the surface shows that the shape of the corrugations changes precisely when passing the vertex of the cone. The reason of this behavior is that $v_1$ the invariant by vertical translation (as opposed to the invariance by central symmetry of the cone and of $\partial_1 f_0$). When $\eta$ decreases toward zero the map $g$ tends towards a piecewise constant map. Each loop in the family $\gamma$ stays at the two points $\cos\theta\,v_1\pm\sin\theta\, n$ for a duration of $\frac{1-\eta}{2}$ each. At the limit, $\gamma$ is a discontinuous map whose image is two points. As a consequence, when $\eta$ is small, the image $f_1(Cyl)$ looks like a piecewise linear surface.

\bibliographystyle{plain}
\bibliography{biblioIC.bib}
\end{document}